\documentclass[12pt]{article}
\usepackage{geometry}
\usepackage{setspace}
\usepackage{fancyhdr}
\usepackage{amsfonts}
\usepackage{amsmath}
\usepackage{verbatim}
\usepackage{listings}
\usepackage{lscape}
\usepackage{graphics}
\usepackage{tabularx}
\usepackage{float}
\usepackage{breqn}
\usepackage{wasysym}
\usepackage{fix-cm}
\usepackage{bbm}
\usepackage{color}
\usepackage{amsthm}
\usepackage{enumitem}
\usepackage{authblk}

\usepackage{appendix}

\usepackage{graphicx}

\usepackage{multirow}
\usepackage{parskip}
\usepackage{amsmath, amsthm, amssymb}
\usepackage{dsfont}
\usepackage{bm}
\usepackage{parskip}

\usepackage{framed,enumitem}
\usepackage{tikz}
\usetikzlibrary{positioning}
\usepackage{accents}
\usepackage{comment}
\usepackage[ruled,vlined]{algorithm2e}
\usepackage{breqn}

\usepackage[round]{natbib}

\newtheorem{lemma}{Lemma}
\newtheorem{mycor}{Corollary}
\newtheorem{myprop}{Proposition}

\theoremstyle{definition}

\makeatletter
\newtheorem*{rep@theorem}{\rep@title}
\newcommand{\newreptheorem}[2]{%
\newenvironment{rep#1}[1]{%
 \def\rep@title{#2 \ref{##1}}%
 \begin{rep@theorem}}%
 {\end{rep@theorem}}}
\makeatother

\DeclareMathOperator*{\argmax}{argmax}

\newcommand{\abs}[1]{\left| #1 \right|}

\def\T{{ \mathrm{\scriptscriptstyle T} }}
\def\pen{{ \mathrm{pen} }}
\def\Pen{{ \mathrm{Pen} }}
\def\tr{{ \mathrm{tr} }}

\newreptheorem{theorem}{Theorem}
\newreptheorem{lemma}{Lemma}
\newreptheorem{mycor}{Corollary}
\newreptheorem{myprop}{Proposition}

\title{Existence of penalised likelihood estimates and posterior propriety of separable prior distributions for Gaussian precision matrices}
\author{Jack Storror Carter}
\affil{Dept. of Economics and Business, Universitat Pompeu Fabra, Spain}
\affil{Data Science Center, Barcelona School of Economics, Spain}
\date{}

\begin{document}

\maketitle

\begin{abstract}
Penalised likelihoods are often used for sparse estimation of a Gaussian precision matrix. In high dimensional settings where the matrix dimension is larger than the sample size, the sample covariance matrix $S$ is not of full rank and the maximum likelihood estimate of the precision matrix does not exist. An additional advantage of some penalised likelihood estimates, for example the graphical lasso, is that it can exist even in such high dimensional settings. This paper gives a thorough analysis of the existence of penalised likelihood estimates for positive semidefinite $S$. Specific tail conditions are provided on the diagonal and off-diagonal penalty functions that ensure existence of the estimate. This is also extended to the Bayesian setting where conditions on separable prior distributions are provided that ensure the resulting posterior distribution is proper.
\end{abstract}

This paper considers the existence of penalised likelihood estimates for a $p$-dimensional Gaussian precision matrix $\Theta = (\theta_{ij})$
\begin{equation}\label{eq:penlike}
    \hat{\Theta}_\Pen = \argmax_{\Theta \succ 0}\, \log\det\Theta - \tr(S \Theta) - \Pen(\Theta)
\end{equation}
where $S$ is a sample covariance matrix, $\Theta \succ 0$ refers to the space of positive definite (PD) matrices, $\tr(A)$ is the trace of the matrix $A$ and $\Pen$ is some penalty function.  We focus on separable penalty functions of the form\footnote{Here the penalty is written only on the upper triangular entries of $\Theta$ due to symmetry. However, many popular penalty functions are written over all entries of $\Theta$ as $\Pen(\Theta) = \sum_{i,j} \tilde{\pen}_{ij}(\theta_{ij})$. Later results can easily be adapted to this form by taking $\pen_{ij} = \tilde{\pen}_{ij} + \tilde{\pen}_{ji}$ for $i < j$.}
\begin{equation}\label{eq:penalty}
   \Pen(\Theta) = \sum_{i \leq j} \pen_{ij}(\theta_{ij}) 
\end{equation}
The diagonal penalty functions $\pen_{ii}$ are defined over the positive real numbers $(0,\infty)$, since this is a requirement for PDness of $\Theta$, and are assumed to be non-negative, non-decreasing and continuous. The off-diagonal penalties $\pen_{ij}$, $i < j$ are defined over all real numbers $(-\infty,\infty)$ and are assumed to be non-negative, symmetric around $0$, non-decreasing in $\abs{\theta_{ij}}$ and continuous. We also assume that there is a common penalty function on all diagonals $\pen_{ii} = \pen_{D},\, i=1,\dots,p$ and on all off-diagonals $\pen_{ij} = \pen,\, i < j$. These are each common assumptions for penalty functions and help with notational simplicity, however it will be shown that the results of this paper can generalise beyond most of these assumptions.

Such penalised likelihood estimates are commonly used for sparse estimation of $\Theta$. The most well known is the graphical lasso (glasso) \citep{Yuan2007,Banerjee2008,Friedman2008} which uses an $l_1$ penalty. Non-convex penalties, such as the SCAD penalty \citep{fan2001variable,fan2009network} and MCP \citep{zhang2010nearly}, were proposed to reduce bias in the estimates of large entries of $\Theta$ in comparison to the glasso. Other penalty functions have been proposed as a continuous approximation of the $l_0$ penalty such as the seamless $l_0$ \citep{dicker2013variable} and ATAN \citep{wang2016variable} penalties. See \cite{williams2020beyond} for the application of these penalties to Gaussian graphical models. 

In the Bayesian framework, $\hat{\Theta}_\Pen$ is equivalent to the maximum a-posteriori (MAP) estimate of $\Theta$ under the prior distribution $p(\Theta) \propto \exp\left(-\frac{n}{2}\Pen(\Theta)\right)$ which has the form
$$p(\Theta) \propto \prod_{i=1}^{p} \pi_D(\theta_{ii}) \prod_{i<j}\pi(\theta_{ij}) \mathbb{I}(\Theta \succ 0)$$
Such prior distributions include the Bayesian glasso \citep{Wang2012}, the graphical horseshoe \citep{li2019graphical} and spike and slab priors \citep{wang2015scaling,Gan2018,jewson2024graphical,sulem2025bayesian}.

When the sample covariance matrix $S$ is PD, the maximum likelihood estimate exists and is equal to $S^{-1}$. It naturally follows that $\hat{\Theta}_\Pen$ exists when $\Pen$ is continuous and non-negative. When $S$ is positive semidefinite (PSD), but not PD (which will be referred to as \textit{only} PSD), the MLE does not exist. However, the addition of a sufficiently strong penalty function can be enough for $\hat{\Theta}_\Pen$ to still exist. This is most well known with the glasso which exists for all PSD $S$ when the diagonal penalty is included \citep[Theorem 1]{Banerjee2008} and for all PSD $S$ with strictly positive diagonal when the diagonal penalty is not included \citep[Theorem 8.7]{lauritzen2022}. However, as will be shown in this paper, many other penalty functions can also achieve such existence properties.

In this paper we will conduct a thorough analysis on the existence of $\hat{\Theta}_\Pen$ for PSD $S$. We provide sufficient tail conditions on $\pen_D$ and $\pen$ that ensure existence of $\hat{\Theta}_\Pen$ for all PSD $S$, as well as conditions that ensure existence with probability 1 when $S$ is a Gaussian sample covariance matrix. We also provide sufficient conditions for $\hat{\Theta}_\Pen$ to not exist for any only PSD $S$. These results naturally extend to the existence of the MAP estimate under the analogous prior. However, a full Bayesian analysis typically requires the stronger condition that the posterior density $p\left(\,\Theta \mid S\,\right)$ is proper, i.e. has finite integral. As such, conditions on $\pi_D$ and $\pi$ will be provided that ensure the posterior $p\left(\,\Theta \mid S\,\right)$ is proper for all PD $S$ and for all PSD $S$.
%and with probability 1 when $S$ is a Gaussian sample covariance matrix.

A summary of the results of the paper and a discussion are in Section \ref{sec:summary}. Background and notation is in Section \ref{sec:bg} and the maximum likelihood estimate is considered in Section \ref{sec:MLE}. In Sections \ref{sec:diag} and \ref{sec:offdiag}, the diagonal penalty $\pen_D$ and the off-diagonal penalty $\pen$ are considered separately while Section \ref{sec:Combined} provides results on joint penalties. Section \ref{sec:gen} gives a lengthy discussion and additional results about how to generalise beyond the above assumptions on the form of $\Pen$. This includes discussion on an alternative form of penalty function which instead penalises partial correlations. Finally, Section \ref{sec:Bayes} considers the Bayesian setting with conditions for posterior propriety.

\section{Summary of results and discussion}\label{sec:summary}

Existence of the penalised likelihood estimate for only PSD $S$ depends on the limiting behaviours of the penalty functions $\pen_D$ and $\pen$. The results in this paper first investigate the required limiting behaviours of $\pen_D$ and $\pen$ individually before giving a joint result.

\subsection{Diagonal penalty}

First considering the diagonal penalty $\pen_D$ combined with any off-diagonal penalty $\pen \geq 0$, $\hat{\Theta}_\Pen$ is guaranteed to exist for any PSD $S$ if (Proposition \ref{prop:diag})
$$\lim_{x \to \infty} \pen_D(x) - \log x = \infty$$
If instead $S$ is a Gaussian sample covariance matrix with $m_0$ eigenvalues equal to 0, then $\hat{\Theta}_\Pen$ exists with probability 1 when (Corollary \ref{cor:DiagProb1})
$$\liminf_{x \to \infty} \frac{\pen_D(x)}{\log x} > m_0 / p$$ 
recalling that $p$ is the matrix dimension.

On the other hand, if $\pen = 0$ and 
$$\liminf_{x \to \infty} p \, \pen_D(x) - \log x = -\infty$$ then $\hat{\Theta}_\Pen$ does not exist for any only PSD $S$ (Proposition \ref{prop:diagnonexist}).

\subsection{Off-diagonal penalty}

When only the off-diagonals are penalised, i.e. $\pen_D=0$, then $\hat{\Theta}_\Pen$ does not exist whenever $S$ has a diagonal entry equal to $0$. This holds no matter the choice of $\pen$. For $\hat{\Theta}_\Pen$ to exist, it is therefore required that $S$ has strictly positive diagonal. Luckily this holds with probability 1 when $S$ is a Gaussian sample covariance matrix.

When $S$ is PSD with strictly positive diagonal, $\hat{\Theta}_\Pen$ is guaranteed to exist when (Corollary \ref{cor:offdiagAll})
$$\liminf_{\abs{x}\to\infty}\frac{\pen(x)}{\log\abs{x}} > p-1$$
If $S$ has $m_0$ eigenvalues equal to $0$ and strictly positive diagonal, then $\hat{\Theta}_\Pen$ exists whenever (Proposition \ref{prop:offdiag}) 
$$\liminf_{\abs{x}\to\infty}\frac{\pen(x)}{\log\abs{x}} > m_0$$
On the other hand, if $\pen_D = 0$ and 
$$\liminf_{\abs{x} \to \infty} q \, \pen(x) - \log\abs{x} = -\infty$$ 
where $q=p(p-1)/2$ is the number of off-diagonals, then $\hat{\Theta}_{\Pen}$ does not exist for any only PSD $S$ (Proposition \ref{prop:offdiagnonexist}).

\subsection{Joint penalties}

When considering both the diagonal penalty $\pen_D$ and the off-diagonal penalty $\pen$, their combination can still ensure existence of $\hat{\Theta}_\Pen$ even when they don't satisfy either of the individual conditions of Corollary \ref{cor:DiagProb1} or Proposition \ref{prop:offdiag}.

If $S$ has $m_0$ eigenvalues equal to 0 and strictly positive diagonal then $\hat{\Theta}_\Pen$ exists whenever (Proposition \ref{prop:joint})
$$\liminf_{x\to\infty}\frac{\pen_D(x)}{\log x} + \liminf_{\abs{x}\to\infty}\frac{\pen(x)}{\log\abs{x}} > m_0$$
On the other hand, if 
$$\liminf_{\abs{x} \to \infty} p \, \pen_D(x) + q \, \pen(x) - \log\abs{x} = -\infty$$ then $\hat{\Theta}_\Pen$ does not exist for any only PSD $S$ (Proposition \ref{prop:NE_joint}).

\subsection{Generalisations}

All the above results extend to unequal penalties where a different penalty $\pen_{ij}$ is applied to each $\theta_{ij}$ (Section \ref{subsec:nonequal}), to non-symmetric $\pen$ (Section \ref{subsec:nonsym}) and to $\pen$ and $\pen_D$ that are not increasing (Section \ref{subsec:notincreasing}). They also trivially extend beyond non-negative penalties to lower bounded penalties. When $\pen_D$ is not lower bounded, the result of Proposition \ref{prop:diag}, can be extended with the condition $\lim_{x \to 0^+} \pen_D(x)-\log x = \infty$ (Proposition \ref{prop:unbounded}).

The results on $\pen_D$ also apply to partial correlation based penalties where the off-diagonal penalty is replaced by a penalty on the partial correlations $-\theta_{ij}/\sqrt{\theta_{ii}\theta_{jj}}$. However, unlike for the standard separable penalty functions, the penalty on the partial correlations alone is not enough to achieve existence when $S$ is PSD (Section \ref{subsec:PC}). 

\subsection{Proper Bayesian posteriors}

The above results on penalty functions also imply the existence of the maximum a-posteriori estimate under the prior $p(\Theta) \propto \exp\left(-\frac{n}{2}\Pen(\Theta)\right) \mathbb{I}(\Theta \succ 0)$. However, slightly stronger tail conditions on $\pi_D$ and $\pi$ are required for the posterior to be proper (i.e. have finite integral) when $S$ is PSD. Even when $S$ is PD, some conditions on $\pi_D$ and $\pi$ are required for posterior propriety, that also depend on the sample size $n$.

If $\pi$ is bounded and $\pi_D$ has $\limsup_{x \to \infty} \pi_D(x) < \infty$ and $\int_{0}^1 x^{(n+p-1)/2} \pi_D(x) \, dx < \infty$ then the posterior is proper for any PD $S$ (Proposition \ref{prop:ProperPD}). If instead $\pi$ is proper and $\pi_D$ has finite $n/2$ moment then the posterior is proper for any PSD $S$ (Proposition \ref{prop:ProperPSD}). 
%The tail condition on $\pi_D$ can be relaxed to $\pi_D(x) \leq C x^{-a}$ with $a > \frac{m_0}{p}\left(\frac{n}{2}+1\right)$ when $S$ is a Gaussian sample covariance matrix with $m_0$ eigenvalues equal to 0. Then the posterior is proper with probability 1 (Proposition \ref{prop:ProperWP1}).

\subsection{Discussion}

The results in this paper give a road map for determining if $\hat{\Theta}_\Pen$ exists for almost any choice of separable $\Pen$. The most useful application of this is in the choice of diagonal penalty $\pen_D$. Since the off-diagonal penalty $\pen$ is the primary tool for inducing sparsity in $\hat{\Theta}_\Pen$, it should be chosen mostly based on statistical and model selection properties. This may be a penalty with weak tail growth, such as the $l_0$ approximating seamless $l_0$ and ATAN penalties, or a bounded penalty such as the SCAD penalty or the MCP. Alone, these penalties on the off-diagonals result in $\hat{\Theta}_\Pen$ not existing whenever $S$ is only PSD and therefore can only be used in real data settings with sufficient sample size $n>p$. This greatly limits their use in high dimensional settings. 

However, the results of this paper show that by combining these penalties with a suitably strong penalty on the diagonals, it can be ensured that $\hat{\Theta}_\Pen$ exists for PSD $S$. Of course, the choice of $\pen_D$ does have an impact on both estimation performance and the sparsity of $\hat{\Theta}_\Pen$. Generally, larger penalisation of the diagonals leads to greater shrinkage of the off-diagonals and greater sparsity. The results in this paper can therefore be used to give some minimal penalisation of the diagonal entries that ensure existence, for example when $\pen_D$ grows logarithmically and exceeds the bound given by Corollary \ref{cor:DiagProb1}.

\section{Background and notation}\label{sec:bg}

The log-likelihood function for a $p \times p$ Gaussian precision matrix $\Theta = (\theta_{ij})$ given a $p \times p$ PSD matrix $S = (s_{ij})$, after removing additive and multiplicative constants, and the corresponding MLE are 
\begin{align*}
    l\left(\,\Theta \mid S\,\right) = \log\det\Theta - \tr(S\Theta) && \hat{\Theta} = \argmax_{\Theta \succ 0}\, l\left(\,\Theta \mid S\,\right)
\end{align*}
where $\tr(A)$ denotes the trace of a matrix $A$ and $\Theta \succ 0$ refers to the set of $p \times p$ PD matrices.

A penalised likelihood subtracts a penalty function from the log-likelihood and estimates $\Theta$ by maximising the resulting function\footnote{Letting $n$ be the sample size, here the multiplicative constant $n/2$ has been removed so that the penalised likelihood function does not depend on $n$. Hence the true penalisation is actually $\frac{n}{2}\Pen(\Theta)$.}
\begin{align*}
    l_{\Pen}\left(\,\Theta \mid S\,\right) = l\left(\,\Theta \mid S\,\right) - \Pen(\Theta) && \hat{\Theta}_{\Pen} = \argmax_{\Theta \succ 0}\, l_{\Pen}\left(\,\Theta \mid S\,\right)
\end{align*}
Recall that in this paper it is assumed that the penalty function is separable with common diagonal and off-diagonal penalty
$$ \Pen(\Theta) = \sum_{i=1}^p \pen_D(\theta_{ii}) + \sum_{i < j} \pen(\theta_{ij})$$
and that $\pen_D$ and $\pen$ are non-negative and continuous and $\pen_D$ is non-decreasing in $\theta_{ii}$ while $\pen$ is symmetric around $0$ and non-decreasing in $\abs{\theta_{ij}}$.

It will be useful to consider the optimisation problem in terms of the eigenvalues and eigenvectors of $\Theta$. Because $\Theta$ is symmetric, it is guaranteed to have an orthonormal basis of eigenvectors. Write the eigenvalues of $\Theta$ as $\lambda = (\lambda_1,\ldots,\lambda_p)$ with corresponding orthonormal eigenvectors $V = (v_1,\ldots,v_p)$. The $j$th entry of the eigenvector $v_k$ is written as $v_{kj}$. For $\Theta$ to be PD it must have strictly positive eigenvalues $\lambda_k > 0$ and its eigenvectors are in the set of orthonormal bases $V \in \mathcal{V}$.

The determinant of a matrix is the product its eigenvalues, and the trace can be written as $\tr(S\Theta) = \sum_{k=1}^{p} \lambda_k v_k^{\top} S v_k$. The log-likelihood function can therefore be rewritten in terms of eigenvalues and eigenvectors as
\begin{align*}
    l\left(\,\Theta \mid S\,\right) = \sum_{k=1}^p \log \lambda_k - \lambda_k v_k^{\top} S v_k
\end{align*}
The penalised likelihood can similarly be rewritten by noting that $\theta_{ij} = \sum_{k=1}^p \lambda_k v_{ki} v_{kj}$.

When $S$ is only PSD it has some eigenvalues equal to $0$ and its null space $\mathcal{N}(S) = \{ v \in \mathbb{R}^p : Sv = 0 \}$ has positive dimension. Let $m_0 = \dim(\mathcal{N}(S))$, which is equal to the number of eigenvalues of $S$ that are equal to $0$. The rank of $S$ is equal to $p - m_0$. We refer to any eigenvector of $\Theta$ in $\mathcal{N}(S)$, $v_k \in \mathcal{N}(S)$ as a null space eigenvector and the corresponding eigenvalue $\lambda_k$ as a null space eigenvalue.

%When $S$ is only PSD, it has some eigenvalues equal to 0, say $\lambda_1,...,\lambda_{m_0} = 0$. In this case, the null space of $S$, denoted by $\mathcal{N}(S) = \{ v \in \mathbb{R}^p : Sv = 0 \}$, is non-empty and of dimension $m_0$. Any $m_0$ orthogonal vectors in $\mathcal{N}(S)$ can be chosen as the eigenvectors of $S$. If $v \in \mathcal{N}(S)$ then it is orthogonal to any vector not in $\mathcal{N}(S)$.

In the optimisation problems, $S$ could be any $p \times p$ PSD matrix. However, it is usually the sample covariance matrix from a $p$-variate Gaussian i.i.d.\ sample $X_1,\ldots,X_n \overset{\mathrm{iid}}{\sim} N_p(\mu,\Theta^{-1})$ with mean vector $\mu$ and covariance matrix $\Theta^{-1}$. When $\mu$ is unknown, the sample covariance matrix is $S = \frac{1}{n}\sum_{i=1}^n(X_i - \bar{X})(X_i - \bar{X})^{\T}$, where $\bar{X} = \frac{1}{n}(X_1 + \cdots + X_n)$. When $n > p$, $S$ is PD with probability 1. However, when $n \leq p$, $S$ is only PSD with $m_0 = p - n + 1$ with probability 1 \citep[Section 8.3]{mathai2022multivariate}.

When $\mu$ is known, the sample covariance matrix is instead $S = \frac{1}{n}\sum_{i=1}^n(X_i - \mu)(X_i - \mu)^{\T}$, which is PD with probability 1 when $n \geq p$, but is only PSD when $n < p$ with $m_0 = p - n$ with probability 1.

\section{Maximum likelihood estimate}\label{sec:MLE}

We begin by considering the existence of the MLE for PD and PSD $S$. While these results hardly need proving, they provide simple examples of the style of proofs that will be used for penalised likelihoods and understanding why the MLE does not exist when $S$ is only PSD helps focus the later proofs.

\begin{myprop}\label{prop:MLE}
    The MLE $\hat{\Theta}$ exists for any PD $S$.
\end{myprop}
\begin{proof}
    Since the likelihood function is continuous, the existence of the MLE follows if $l\left(\,\Theta \mid S\,\right) \to -\infty$ whenever $\Theta$ approaches the boundary of the space of PD matrices or $\lVert \Theta \rVert \to \infty$. PD matrices are characterised by positive eigenvalues $\lambda_1,\ldots,\lambda_p > 0$, so the boundary of the space occurs when any $\lambda_k \to 0$ while $\lVert \Theta \rVert \to \infty$ requires that some $\lambda_k \to \infty$. To prove that $l\left(\,\Theta \mid S\,\right) \to -\infty$ in such cases, we will upper bound the objective function by an expression that holds uniformly over all $V \in \mathcal{V}$. %which is a compact set. 
    We will then show that this upper bound tends to $-\infty$ whenever any (potentially multiple) $\lambda_k \to 0,\infty$.

    Because the log-likelihood is separable in the $\lambda_k$, each $\lambda_k$ can be considered separately. Since $S$ is PD, for any unit vector $v$, $v^\top S v \geq \lambda_\min(S) > 0$, where $\lambda_\min(S)$ is the smallest eigenvalue of $S$.
    %$S$ is PD so it has strictly positive eigenvalues $\lambda_1,\ldots,\lambda_p > 0$. Letting $\lambda_{min} = \min_{j=1,\dots,p} \lambda_j$ and noting that $\sum_{j=1}^p (w_k^\T v_j)^2 = 1$ because $(v_1,\dots,v_p)$ is an orthonormal basis and $\lVert w_k \rVert = 1$, we have $\sum_{j=1}^p \lambda_j (w_k^\T v_j)^2 \geq \lambda_{min}\sum_{j=1}^p (w_k^\T v_j)^2 = \lambda_{min} > 0$. 
    Hence the $\lambda_k$ component of the log-likelihood has upper bound
    $$\log\lambda_k - \lambda_k v_k^\top S v_k \leq \log\lambda_k - \lambda_k \lambda_{min}(S)$$
    and this upper bound tends to $-\infty$ as $\lambda_k \to 0,\infty$.
\end{proof}

Since the log-likelihood $l\left(\,\Theta \mid S\,\right)$ is strictly concave in $\Theta \succ 0$, it follows that, when it exists, the MLE is unique and it is well known that when $S$ is PD the MLE is $\hat{\Theta} = S^{-1}$. Since it is assumed that $\Pen \geq 0$ and is continuous, existence of the penalised likelihood estimate also trivially follows for any PD $S$.

\begin{mycor}\label{cor:PD}
    For any non-negative, continuous penalty function $\Pen$, the penalised likelihood estimate $\hat{\Theta}_{\Pen}$ exists for all PD $S$.
\end{mycor}

On the other hand, when $S$ is only PSD, it is not invertible and the MLE does not exist.

\begin{myprop}
    The MLE does not exist when $S$ is only PSD.
\end{myprop}
\begin{proof}
    To prove non-existence we construct a sequence of $\Theta$ such that $l\left(\,\Theta \mid S\,\right) \to \infty$. Since $S$ is only PSD, $\dim(\mathcal{N}(S)) = m_0 > 0$. Let $v_1 \in \mathcal{N}(S)$ be an eigenvector of $\Theta$. Then $v_1^\top S v_1 = 0$ and so the contribution of $\lambda_1$ to the log-likelihood is $\log \lambda_1 - \lambda_1 v_1^\top S v_1 = \log\lambda_1 \to \infty$ as $\lambda_1 \to \infty$. Then, keeping $\lambda_2,\ldots,\lambda_p > 0$ and $v_1,\dots,v_p$ fixed, as $\lambda_1 \to \infty$, $l\left(\,\Theta \mid S\,\right) \to \infty$ and so the MLE does not exist. 
\end{proof}

This proof highlights that the trace term does not depend on $\lambda_k$ whenever $v_k \in \mathcal{N}(S)$ and it is exactly the divergence of these eigenvalues that cause $l\left(\,\Theta \mid S\,\right) \to \infty$. On the other hand, when $v_k \notin \mathcal{N}(S)$ the trace term still grows linearly with $\lambda_k$. Proving the existence of the penalised likelihood estimate $\hat{\Theta}_\Pen$ would usually require showing that $l_\Pen\left(\,\Theta^{(t)} \mid S\,\right) \to -\infty$ for any sequence of PD matrices with $\Theta^{(t)}$ approaching the boundary of the PD cone or $\lVert \Theta^{(t)} \rVert \to \infty$. However, these observations about the likelihood function restrict the types of sequences that must be considered.

\begin{lemma}\label{lemma}
    Let $S$ be PSD with $\dim(\mathcal{N}(S)) = m_0$ and $\Pen$ be a penalty function satisfying the assumptions given in Section \ref{sec:bg}. If $\l_\Pen\left(\,\Theta^{(t)} \mid S\,\right) \to -\infty$ for any sequence of $\Theta^{(t)}$ with eigenvalues $\lambda_1^{(t)},\dots,\lambda_p^{(t)}$ and eigenvectors $V^{(t)}=(v_1^{(t)},\dots,v_p^{(t)})$ satisfying
    \begin{enumerate}[label=(\roman*)]
        \item\label{item1} $\lambda_\max^{(t)} = \max_{k=1,\dots,p}\lambda_k^{(t)} \to \infty$.
        \item\label{item2} $V^{(t)} \to V$ with $v_1,\dots,v_m \in \mathcal{N}(S)$ and $v_{m+1},\dots,v_p \notin \mathcal{N}(S)$, $m \in \{1,\dots,m_0\}$, i.e. the sequence of eigenvectors converge to a basis with $m \in \{1,\dots,m_0\}$ vectors in $\mathcal{N}(S)$.
        %\item $\inf_t \lambda_\min^{(t)} > 0$ where $\lambda_\min^{(t)} = \min_{k=1,\dots,p} \lambda_k^{(t)}$, i.e. all eigenvalues remain away from 0.
        \item\label{item3} $\lambda_k^{(t)} = o(\lambda_\max^{(t)})$ for $k=m+1,\dots,p$, i.e. we only have to consider sequences where non-null space eigenvalues are of smaller order than $\lambda_\max^{(t)}$. Note that this implies that for sufficiently large $t$, $\lambda_\max^{(t)} = \max_{k=1,\dots,m}\lambda_k^{(t)}$ is one of the null space eigenvalues.
    \end{enumerate}
    Then $\hat{\Theta}_\Pen$ exists.
\end{lemma}
\begin{proof}
    As noted previously, showing that $l_\Pen\left(\,\Theta^{(t)} \mid S\,\right) \to -\infty$ for any sequence $\Theta^{(t)}$ approaching the boundary of the PD cone or with $\lVert \Theta^{(t)} \rVert \to \infty$ is enough to prove the existence of $\hat{\Theta}_\Pen$. Such sequences are exactly characterised by at least one of $\lambda_\min^{(t)} = \min_{k=1,\dots,p}\lambda_k^{(t)} \to 0$ or $\lambda_\max^{(t)} \to \infty$. Proving this lemma therefore requires showing that $l_\Pen\left(\,\Theta^{(t)} \mid S\,\right) \to -\infty$ for any $\Theta^{(t)}$ satisfying \ref{item1}-\ref{item3} implies that $l_\Pen\left(\,\tilde\Theta^{(t)} \mid S\,\right) \to -\infty$ for any $\tilde{\Theta}^{(t)}$ with $\tilde{\lambda}_\min^{(t)} \to 0$ or $\tilde{\lambda}_\max^{(t)} \to \infty$. We do this for each of \ref{item1}-\ref{item3} individually.
    %We show that $\l_\Pen\left(\,\Theta^{(t)} \mid S\,\right) \to -\infty$ for any $\Theta^{(t)}$ satisfying properties \ref{item1}-\ref{item3} implies that $\l_\Pen(\tilde{\Theta}^{(t)} \mid S) \to -\infty$ for general $\tilde{\Theta}^{(t)}$ approaching the boundary of the PD cone or with $\lVert \tilde{\Theta}^{(t)} \rVert \to \infty$. We do this by building up one property at a time.
    \begin{enumerate}[label=(\roman*)]
        \item The proof of Proposition \ref{prop:MLE} already demonstrates that $l\left(\,\Theta^{(t)} \mid S\,\right) \to -\infty$ whenever $\lambda_\min^{(t)} \to 0$ but $\lambda_\max^{(t)} \not\to \infty$. This is regardless of $S$ or $V^{(t)}$. Since $\Pen \geq 0$ it follows that $l_\Pen\left(\,\Theta^{(t)} \mid S\,\right) \to -\infty$ also.

        All remaining sequences that approach the boundary of the PD cone or have $\lVert \Theta^{(t)} \rVert \to \infty$ must have $\lambda_\max^{(t)} \to \infty$. Note that this does not rule out the case where both $\lambda_\min^{(t)} \to 0$ and $\lambda_\max^{(t)} \to \infty$.
        
        \item First note that the case where $m=0$ eigenvectors converge in $\mathcal{N}(S)$ is analogous to Corollary \ref{cor:PD} while the case where $m > m_0$ is impossible because $\dim(\mathcal{N}(S)) = m_0$ and so only $m=\{1,\dots,m_0\}$ must be considered.
        
        Suppose that $l_\Pen\left(\,\Theta^{(t)} \mid S\,\right) \to -\infty$ for any sequence of $\Theta^{(t)}$ satisfying \ref{item1} and \ref{item2}. Consider a general sequence $\tilde{\Theta}^{(t)}$ with $\tilde{\lambda}_\max^{(t)} \to \infty$ and suppose that $l_\Pen\left(\,\tilde{\Theta}^{(t)} \mid S\,\right) \not\to -\infty$. Then there is a subsequence $\tilde{\Theta}^{(t_i)}$ such that $l_\Pen\left(\,\tilde\Theta^{(t)} \mid S\,\right) \geq C$ for all $i$ and some $C \in \mathbb{R}$. Since the eigenvectors $\tilde{V}^{(t_i)}$ lie in the compact set of orthonormal bases, there is a further subsequence with $V^{(j)} = \tilde{V}^{(t_{i_j})} \to V^*$ converging. The resulting subsequence $\Theta^{(j)}$ then satisfies property \ref{item2} and so $l_\Pen\left(\,\Theta^{(j)} \mid S\,\right) \to -\infty$ which gives a contradiction. Hence $l_\Pen\left(\,\tilde\Theta^{(t)} \mid S\,\right) \to -\infty$.
        
        %\item Now assume that $l_\Pen(\Theta^{(t)} \mid S) \to -\infty$ for any sequence of $\Theta^{(t)}$ satisfying properties \ref{item1},\ref{item2} and with $\lambda_\max^{(t)} \to \infty$. Consider a sequence $\tilde{\Theta}^{(t)}$ satisfying property \ref{item1} and with $\tilde{\lambda}_\min^{(t)} \to 0$. If $\tilde{\lambda}_\max^{(t)} \not\to \infty$, then, as seen in Proposition \ref{prop:MLE}, the log likelihood tends to $-\infty$ as $\tilde{\lambda}_\min^{(t)} \to 0$. Because $\Pen \geq 0$, this implies that $l_\Pen(\tilde\Theta^{(t)} \mid S) \to -\infty$. If $\tilde{\lambda}_\max^{(t)} \to \infty$ then...
        
        %For a sequence of $\tilde{\lambda}_1^{(t)},\dots,\tilde{\lambda}_p^{(t)}$ with $\min_{k=1,\dots,p} \tilde{\lambda}_k^{(t)} \to 0$, consider an alternative sequence $\lambda_k^{(t)} = \max\{\epsilon,\tilde{\lambda}^{(t)}\}$ for $\epsilon > 0$. As seem in the proof of Proposition \ref{prop:MLE}, the $\tilde{\lambda}_k$ component of the log-likelihood tends to $-\infty$ as $\tilde{\lambda}_k \to 0$ and, for sufficiently small $\epsilon > 0$, is increasing on $\tilde{\lambda}_k \in (0,\epsilon)$. Hence $l\left(\tilde{\Theta}^{(t)} \mid S\right) \leq l\left(\Theta^{(t)} \mid S\right)$. Additionally, the difference $\lVert \Theta^{(t)} - \tilde{\Theta}^{(t)} \rVert_2 \leq \epsilon$ is uniformly bounded and so clipping small eigenvalues does not affect the entry-wise divergence responsible for the penalty asymptotics. This implies that if $\l_\Pen(\Theta^{(t)} \mid S) \to -\infty$ then $\l_\Pen(\tilde{\Theta}^{(t)} \mid S) \to -\infty$ also.

        \item Now suppose that $l_\Pen\left(\,\Theta^{(t)} \mid S\,\right) \to -\infty$ for any sequence of $\Theta^{(t)}$ satisfying \ref{item1}-\ref{item3}. Consider a sequence $\tilde{\Theta}^{(t)}$ satisfying \ref{item1} and \ref{item2} but with $\lambda_k^{(t)} \neq o\left(\lambda_\max^{(t)}\right)$ for some $k=m+1,\dots,p$. This implies that $\tilde{\lambda}_k^{(t)} \geq \epsilon \, \tilde{\lambda}_\max^{(t)}$ infinitely often for some $\epsilon > 0$. Consider any subsequence $\tilde{\Theta}^{(t_i)}$ along which $\tilde{\lambda}_k^{(t_i)} \geq \epsilon \, \tilde{\lambda}_\max^{(t_i)}$ for some $k=m+1,\dots,p$ and some $\epsilon > 0$. Because $v_k \notin \mathcal{N}(S)$, eventually $v_k^{(t)\top} S v_k^{(t)} \geq c_k > 0$ and so $\tr(S\Theta^{(t)}) \geq c_k \lambda_k^{(t)}$. Using this, replacing all eigenvalues with $\lambda_\max^{(t)}$ in the log determinant and removing the penalty term,
        \begin{align*}
        l_\Pen\left(\,\Theta^{(t_i)} \mid S\,\right) &\leq p \log \lambda_\max^{(t_i)} - c_k \lambda_k^{(t_i)} \\
        & \leq p \log \lambda_\max^{(t_i)} - c_k \, \epsilon \, \lambda_\max^{(t_i)}
        \end{align*}
        This upper bound tends to $-\infty$ as $\lambda_\max^{(t_i)} \to \infty$ and so $l_\Pen\left(\,\Theta^{(t_i)} \mid S\,\right) \to -\infty$.

        The complement of these subsequences has $\tilde{\lambda}_k^{(t)} < \epsilon \, \tilde{\lambda}_\max^{(t)}$ for all $k=m+1,\dots,p$ and all $\epsilon > 0$. If these form a subsequence $\Theta^{(j)} = \tilde{\Theta}^{(t_j)}$, then it satisfies \ref{item1}-\ref{item3} and so $l_\Pen\left(\,\Theta^{(j)} \mid S\,\right) \to -\infty$.

        It follows that $l_\Pen\left(\,\tilde\Theta^{(t)} \mid S\,\right) \to -\infty$ along the full sequence.
    \end{enumerate}
\end{proof}
This lemma shows that when proving the existence of $\hat\Theta_\Pen$ it is enough to only consider sequences that satisfy \ref{item1}-\ref{item3}. In particular, we only need to consider sequences in which the eigenvectors converge to a basis with at least one null-space vector, and sequences of eigenvalues with at least one null space eigenvalue tending to $\infty$.

\section{Diagonal penalty}\label{sec:diag}

While the MLE does not exist for only PSD $S$, the addition of a sufficiently strong penalty on the diagonal is enough to ensure existence. In this section we remove the off-diagonal penalty $\pen = 0$, but all existence results remain valid for any $\pen \geq 0$.

\begin{myprop}\label{prop:diag}
    If $\lim_{x \to \infty} \pen_D(x) - \log x = \infty$, then the penalised likelihood estimate $\hat{\Theta}_{\Pen}$ exists for any PSD $S$.
\end{myprop}
\begin{proof}
    Since $\tr(S\Theta) \geq 0$ for PSD $S$ and PD $\Theta$, the penalised likelihood is upper bounded by removing the trace term (this corresponds to a worst case scenario when $S$ is a zero matrix and all $v_k \in \mathcal{N}(S)$). Hadamard's inequality also says that the determinant of $\Theta$ is upper bounded by the product of its diagonal entries and therefore $\sum_{i=1}^p \log \lambda_i  \leq \sum_{i=1}^p \log\theta_{ii}$. These give the upper bound
    \begin{align*}
        l_{\Pen}\left(\,\Theta \mid S \,\right) 
        %&\leq \sum_{i=1}^p \log \lambda_i - \pen_D\left( \theta_{ii} \right) \\
        &\leq \sum_{i=1}^p \log(\theta_{ii}) - \pen_D(\theta_{ii})
    \end{align*}
    Since $v_1,\ldots,v_p$ are unit-vectors, for each $k=1,\ldots,p$ we have $\max_{j=1,\dots,p} (v_{kj})^2 \geq 1/p$. By choosing $i = i(k) \in \argmax_{j=1,\dots,p} (v_{kj})^2$, we have $\theta_{ii} \geq \lambda_k (v_{ki})^2 \geq \lambda_k / p$. In particular, there is an $i$ such that $\theta_{ii} \geq \lambda_\max/p$ and so as $\lambda_\max \to \infty$, $\theta_{ii} \to \infty$ also. It follows that as $\lambda_\max \to \infty$ we have $\sum_{i=1}^p \log(\theta_{ii}) - \pen_D(\theta_{ii}) \to -\infty$ and therefore $l_{\Pen}\left(\,\Theta \mid S \,\right) \to -\infty$.
\end{proof}

On the other hand, under a weaker diagonal penalty with $\liminf_{x \to \infty} \pen_D(x) - \log x = -\infty$, the solution does not exist for some PSD $S$. This is most easily seen by considering $S = 0$ and $\Theta$ with standard basis eigenvectors in which case the penalised likelihood tends to $\infty$ as any $\lambda_k \to \infty$. However, it does not guarantee that the solution does not exist for all only PSD $S$. In particular, if multiple $\theta_{ii} \to \infty$ as $\lambda_k \to \infty$, then these multiple penalty terms might be enough to overpower the $\log$ term. But an even slower growth condition on $\pen_D$ does ensure this.

\begin{myprop}\label{prop:diagnonexist}
    If $\pen = 0$ and $\liminf_{x \to \infty} p \, \pen_D(x) - \log x = -\infty$, then $\hat{\Theta}_{\Pen}$ does not exist for any only PSD $S$.
\end{myprop}
\begin{proof}
    Let $S$ be only PSD so that $\dim(\mathcal{N}(S)) = m_0 > 0$. Choose a unit vector $v_1 \in \mathcal{N}(S)$ and let $\Theta^{(t)} = I_p + tv_1v_1^\top$ for $t > 0$. Then $\Theta^{(t)}\succ 0$ for all $t > 0$, $\det(\Theta^{(t)})=1+t$, and since $v \in N(S)$ we have
    $\tr(S\Theta^{(t)})=\tr(S I_p)+t\,v^\top S v=\tr(S)$, which is constant in $t$. The diagonal entries of $\Theta^{(t)}$ are $\theta_{ii}^{(t)} = 1+tv_{1i}^2 \leq 1+t$. Since $\pen_D$ is increasing, the penalty term is upper bounded by
    \begin{align*}
        \sum_{i=1}^p \pen_D(\theta_{ii}^{(t)}) 
        %&\leq \sum_{i=1}^p \pen_D(1+t) \\
        &\leq p \, \pen_D(1+t)
    \end{align*}
    The penalised likelihood function is therefore lower bounded by 
    $$ l_{\Pen} \left(\,\Theta^{(t)} \mid S \,\right) \geq \log(1+t) - \tr(S) - p \, \pen_D(1+t)$$
    The condition $\liminf_{x\to\infty} p\,\pen_D(x) - \log x = -\infty$ ensures that there is an increasing sequence of $t$ such that this lower bound tends to $\infty$ and therefore along this sequence, $l_{\Pen} \left(\,\Theta^{(t)} \mid S \,\right) \to \infty$.
    %The condition $\liminf_{x \to \infty} p \, \pen_D(x) - \log x = -\infty$ implies that for every $M > 0$ there is an $x$ such that $\log x - p \, \pen_D(x) \geq M$. Choosing a sequence $x_m$ such that $\log x_m - p \, \pen_D(x_m) \geq m$ and letting $t_m = x_m -1$ gives
    %\begin{align*}
    %    l_{\Pen} \left( \Theta^{(t_m)} \mid S \right) &\geq \log x_m - p \, \pen_D(x_m) - \tr(S) \\
    %    &\geq m - \tr(S)
    %\end{align*}
    %and this lower bound tends to $\infty$ as $m \to \infty$. Hence $l_{\Pen} \left( \Theta \mid S \right)$ is unbounded above over $\Theta \succ 0$ and so the maximiser $\hat{\Theta}_{\Pen}$ does not exist.
\end{proof}

We have now seen that when $\lim_{x \to \infty} \pen_D(x) - \log x = \infty$, the solution exists for every PSD $S$, while if $\liminf_{x \to \infty} p \, \pen_D(x) - \log x = -\infty$ the solution does not exist for any only PSD $S$. Both of these conditions and results are independent of the rank of $S$ and do not require that $S$ be a Gaussian sample covariance matrix. This begs the question of what happens when $\pen_D$ satisfies neither of these conditions. Is there some other condition on $\pen_D$ that guarantees existence for some useful subset of PSD $S$? 

One specific such class of penalty functions $\pen_D(x) = c \, \log x$ was investigated by \cite{carter2025existence} where it was shown that the solution exists with probability $1$ when $S$ is a Gaussian sample covariance matrix for certain values of $c$ depending on the rank of $S$. We restate the result of \cite{carter2025existence}, Theorem 1.1 in the context of this paper.

\begin{myprop}\label{prop:otherpaper}
    Let $\pen_D(x) = c \, \log x$ and $S$ be a Gaussian sample covariance matrix with $m_0$ eigenvalues equal to $0$. Then $\hat{\Theta}_\Pen$ exists with probability 1 for any $c \in \left(\frac{m_0}{p},1\right)$.
\end{myprop}
This result relies on the fact that orthogonal vectors in $\mathcal{N}(S)$ can have limited sparsity with probability 1 - this will be explored further in Section \ref{sec:Bayes}. The penalty parameter $c$ must exceed $m_0/p$ so that there is sufficient penalisation in each direction of $\mathcal{N}(S)$. Also note that the penalty $\pen_D(x) = c \log x$ is not lower bounded and so does not meet the previous assumptions - this will be discussed further in Section \ref{subsec:nonneg} and is the reason why $c < 1$ is required. 

The proof of Proposition \ref{prop:otherpaper} in \citet{carter2025existence} follows the same format as the proofs in this paper, showing $l_\Pen\left(\,\Theta \mid S \,\right) \to -\infty$ whenever $\Theta$ approaches the boundary of the PD cone or $\lVert \Theta \rVert \to \infty$. This leads to a simple generalisation of the result.

\begin{mycor}\label{cor:DiagProb1}
    Let $S$ be a Gaussian sample covariance matrix with $m_0$ eigenvalues equal to $0$. If
    $$\liminf_{x \to \infty} \frac{\pen_D(x)}{\log x} > \frac{m_0}{p}$$
    then $\hat{\Theta}_\Pen$ exists with probability 1
\end{mycor}
\begin{proof}
    Suppose 
    $$\liminf_{x \to \infty} \frac{\pen_D(x)}{\log x} = c_0 > \frac{m_0}{p}$$
    and let $c \in (m_0/p,\min\{c_0,1\})$. Then there exists $x_0$ such that $\pen_D(x) \geq c \, \log x$ for all $x > x_0$. It can then be globally lower bounded by $\pen_D(x) \geq c \, \log x - a$ for all $x>0$ where $a = \sup_{0 < x \leq x_0}(c \log x - \pen_D(x))$. Note that $a$ is a finite constant because $x_0$ is fixed, the function is continuous on $(0,x_0]$ and $\pen_D \geq 0$ so $c \log x - \pen_D(x) \to -\infty$ as $x \to 0$. Then we have $\sum_{i=1}^p \pen_D(\theta_{ii}) \geq \sum_{i=1}^p c \, \log \theta_{ii} - pa$. It follows that
    $$ l_\Pen\left(\,\Theta \mid S \,\right) \leq l_{\Pen'}\left(\,\Theta \mid S \,\right) + pa$$
    where $\Pen'$ has $\pen_D'(x) = c\log x$ and $\pen' = 0$ so is a penalty of the form of Proposition \ref{prop:otherpaper} and $pa$ is a constant. The result follows.
\end{proof}

This section has covered almost all tail conditions of $\pen_D$. If $\liminf_{x \to \infty} p \, \pen_D(x) - \log x = -\infty$ then the solution does not exist for any only PSD $S$ while if $\liminf_{x \to \infty} \pen_D(x) - \log x = -\infty$ it does not exist for some PSD $S$. If $\lim_{x \to \infty} \pen_D(x) - \log x = \infty$ then the solution exists for all PSD $S$ while if $\liminf_{x \to \infty} \pen_D(x)/\log x > m_0/p$ then the solution exists with probability 1 for $S$ of rank $p-m_0$. This extends to when $\liminf_{x \to \infty} \pen_D(x)/\log x > 1/p$ when the solution exists with probability 1 for $S$ of rank $p-1$. The only tail conditions remaining are the boundary conditions
\begin{itemize}
    \item $\liminf_{x \to \infty} \pen_D(x)/\log x = 1/p$ exactly (because $\liminf_{x \to \infty} \pen_D(x)/\log x < 1/p$ implies $\liminf_{x \to \infty} p \, \pen_D(x) - \log x = -\infty$). Our results show that in this case there are some PSD $S$ for which $\hat{\Theta}_\Pen$ does not exist. However, it is possible that the solution still exists for some choices of only PSD $S$ and $\pen_D$.
    \item $\liminf_{x \to \infty} \pen_D(x) - \log x \in (-\infty,\infty)$. Our results show that $\hat{\Theta}_\Pen$ exists with probability 1 in this case. However, it is undetermined if the solution exists for all PSD $S$. This may be the case for some such $\pen_D$, for example when $\liminf_{x \to \infty} \pen_D(x) - \log x \gg 0$, but not for others, when $\liminf_{x \to \infty} \pen_D(x) - \log x \ll 0$.
\end{itemize}

\section{Off-diagonal penalty}\label{sec:offdiag}

The previous section focused on the diagonal penalty. In this section the diagonal penalty is instead removed $\pen_D = 0$ so that the penalty is only on the off-diagonals through $\pen$.

%Working with only the off-diagonal penalty introduces some new dyniamics. First, Lemma \ref{lemma} can be updated to also exclude any sequence with $\lambda_\min \to 0$.

%\begin{lemma}\label{lemma2}
%    If $\pen_D = 0$ then Lemma \ref{lemma} can be extended to sequences $\Theta^{(t)}$ satisfying \ref{item1}-\ref{item3} and
%    \begin{enumerate}[label=(\roman*)]
%    \setcounter{enumi}{3}
%        \item\label{item4} $\lambda_\min^{(t)} \geq \epsilon > 0$, i.e. all eigenvalues are bounded away from 0.
%    \end{enumerate}
%\end{lemma}
%\begin{proof}
%    Let $\tilde{\Theta}^{(t)}$ be a sequence satisfying \ref{item1}-\ref{item3} such that for any $\epsilon > 0$ there is a $t$ such that $\tilde{\lambda}_\min^{(t)} < \epsilon$. Consider $\Theta^{(t)} = \tilde{\Theta}^{(t)} + \epsilon I$ for some $\epsilon > 0$. The eigenvalues of $\Theta^{(t)}$ are $\lambda_k^{(t)} = \tilde{\lambda}_k^{(t)} + \epsilon$ while the eigenvectors are $V^{(t)} = \tilde{V}^{(t)}$ and so $\Theta^{(t)}$ satisfies \ref{item1}-\ref{item3} and \ref{item4}. Since $\Theta^{(t)}$ and $\tilde{\Theta}^{(t)}$ have the same off-diagonal entries, $\Pen(\Theta^{(t)}) = \Pen(\tilde{\Theta}^{(t)})$. We also have $\log \det \tilde{\Theta}^{(t)} < \log \det \Theta^{(t)}$ and $\tr(S \Theta^{(t)}) = \tr(S \tilde{\Theta}^{(t)}) + \epsilon \tr(S)$ and so
%    $$ l_\Pen(\tilde{\Theta}^{(t)} \mid S) < l_\Pen(\Theta^{(t)} \mid S) + \epsilon \tr(S)$$
%    Hence $l_\Pen(\Theta^{(t)} \mid S) \to -\infty$ implies that $l_\Pen(\tilde{\Theta}^{(t)} \mid S) \to -\infty$ also.
%\end{proof}

Using only the off-diagonal penalty introduces an additional complication. An off-diagonal entry of $\Theta$ in terms of the eigenvalues and eigenvectors is $\theta_{ij} = \sum_{k=1}^p \lambda_k v_{ki} v_{kj}$. By taking $v_k$ to be a standard basis vector with only one non-zero entry, $v_{ki} v_{kj} = 0$ for all $i \neq j$, none of the off-diagonals depend on the value of $\lambda_k$ and so the penalty function remains finite as $\lambda_k \to \infty$. If $v_k \in \mathcal{N}(S)$, then the trace term also doesn't depend on $\lambda_k$ and so the penalised likelihood is unbounded as $\lambda_k \to \infty$. It is therefore necessary that standard basis vectors do not appear in $\mathcal{N}(S)$ for existence to be possible with only off-diagonal penalisation. Fortunately, the null space of a PSD matrix can only contain a standard basis vector in a very specific circumstance.

\begin{lemma}\label{lem}
    Let $S$ be a PSD matrix. Then $\mathcal{N}(S)$ contains a standard basis vector if and only if $S$ has a diagonal entry equal to 0.
\end{lemma}
\begin{proof}
    First suppose, w.l.o.g., that $e_1 = (1,0,\dots,0)^\top \in \mathcal{N}(S)$. Then $S e_1 = 0$ and therefore $S_{11}=0$.
    
    Now suppose, w.l.o.g., that $S_{11}=0$. Then for $S$ to be PSD it is required that $S_{1i}=S_{i1}=0$ for all $i=1,...,p$. It follows that $S e_1 = 0$ and so $e_1 \in \mathcal{N}(S)$.
\end{proof}

Therefore, when $S$ has strictly positive diagonal entries, $\mathcal{N}(S)$ does not contain any standard basis vectors. It follows that either the trace term or the penalty term is a function of each eigenvalue of $\Theta$.

While this implies that at least one $\theta_{ij}$ depends on the null space eigenvalues, it does not imply that $\abs{\theta_{ij}} \to \infty$ as the eigenvalues grow. By allowing the eigenvalues to grow at different rates, it is possible to keep any $\abs{\theta_{ij}}$ constant. However, it is not possible to do this for all $\theta_{ij}$.

\begin{lemma}\label{lem2}
    Let $S$ be PSD with strictly positive diagonal entries and $\dim(\mathcal{N}(S)) = m_0$. Let $\Theta^{(t)}$ be a sequence of PD matrices as in Lemma \ref{lemma} and $\theta_{\max}^{(t)} = \theta_{i_tj_t}^{(t)}$ where $(i_t,j_t) \in \argmax_{i < j}\abs{\theta_{ij}^{(t)}}$ be the off-diagonal entry of $\Theta^{(t)}$ with largest absolute value.
    Then there exists $\delta>0$ and $t_0$ such that $$\abs{\theta_{\max}^{(t)}} \geq \delta \, \lambda_{\max}^{(t)}$$ for all $t \geq t_0$.
\end{lemma}
\begin{proof}
    Decompose $\Theta^{(t)}$ into the null space and non-null space parts 
    $$\Theta^{(t)} = \sum_{k=1}^{m}\lambda_k^{(t)} v_k^{(t)} v_k^{(t)\top} + \sum_{k=m+1}^p \lambda_k^{(t)} v_k^{(t)} v_k^{(t)\top} = A^{(t)} + C^{(t)}$$
    By letting $\alpha_k^{(t)} = \lambda_k^{(t)}/\lambda_{\max}^{(t)} \in [0,1]$ for $k=1,\dots,m$, we decompose the matrix $A^{(t)}$ further as
    \begin{align*}
    A^{(t)} = \lambda^{(t)}_{\max} B^{(t)} && B^{(t)} = \sum_{k=1}^{m} \alpha_k^{(t)} v_k^{(t)} v_k^{(t)\top}
    \end{align*}
    The off-diagonals of $\Theta$ are then written as $\theta_{ij}^{(t)} = \lambda_{\max}^{(t)} b_{ij}^{(t)} + c_{ij}^{(t)}$. Let $b_{\max}^{(t)} = \max_{i<j} \abs{b_{ij}^{(t)}}$ and $c_{\max}^{(t)} = \max_{i<j}\abs{c_{ij}^{(t)}} < \infty$. This leads to the lower bound
    $$\abs{\theta_{\max}^{(t)}} \geq \lambda_{\max}^{(t)} b_{\max}^{(t)} - c_{\max}^{(t)}$$
    Note that $C^{(t)}$ has eigenvalues $\lambda_{m+1}^{(t)},\dots,\lambda_p^{(t)}$, while the remaining eigenvalues are repeated $0$. Hence $c_{\max}^{(t)} \leq \lVert C^{(t)} \rVert_2 = \max_{k=m+1,\dots,p}\lambda_k^{(t)} = o\left(\lambda_\max^{(t)}\right)$, with the final equality coming from property \ref{item3} of Lemma \ref{lemma}.

    Next we show that $b_\max^{(t)}$ is eventually bounded away from $0$. Suppose it is not so that there is a subsequence such that $b_\max^{(t_i)} \to 0$. Passing to a further subsequence such that $\alpha_k^{(t_{i_j})} \to \alpha_k \in [0,1]$ for each $k=1,\dots,m$ and because $v_k^{(t)} \to v_k$, we have $B^{(t_{i_j})} \to B = \sum_{k=1}^{m} \alpha_k v_k v_k^{\top}$. Note that since $v_1,\dots,v_m \in \mathcal{N}(S)$, we have $\mathrm{range}(B) \subseteq \mathcal{N}(S)$. Also, since $\max_{k=m+1,\dots,p}\lambda_k^{(t)} = o\left(\lambda_\max^{(t)}\right)$, eventually $\lambda_\max^{(t)} = \max_{k=1,\dots,m} \lambda_k^{(t)}$ and so $\max_{k=1,\dots,m} \alpha_k^{(t)} = 1$. Hence $B^{(t)}$ eventually has maximum eigenvalue equal to 1 and $\lVert B^{(t)} \rVert_2 = 1$. Hence $\lVert B \rVert_2 = 1$ and so $B \neq 0$. Since $b_\max^{(t_i)} \to 0$, all off-diagonals of $B$ must be equal to 0 and, by $B \neq 0$, $B$ must have at least one positive diagonal entry. It follows that a basis vector $e_r \in \mathrm{range}(B) \subseteq \mathcal{N}(S)$. By Lemma \ref{lem} this contradicts $S$ having strictly positive diagonal. Hence $\liminf_t b_\max^{(t)} > 0$ and so there is $\delta_0 > 0$ and $t_1$ such that $b_{\max}^{(t)} \geq \delta_0 > 0$ for all $t \geq t_1$.

    We therefore have 
    $$\abs{\theta_{\max}^{(t)}} \geq \lambda_{\max}^{(t)} \delta_0 - c_\max^{(t)}$$
    for all $t \geq t_1$. Since $c_\max^{(t)} = o\left(\lambda_{\max}^{(t)}\right)$, there exists $t_2$ such that $\lambda_{\max}^{(t)} \delta_0 \geq 2c_{\max}^{(t)}$ for all $t \geq t_2$. Letting $t_0 = \max\{ t_1,t_2 \}$, we have for all $t \geq t_0$,
    \begin{align*}
        \abs{\theta_{\max}^{(t)}} &\geq \lambda_{\max}^{(t)} \delta_0 - c_{\max}^{(t)}\\
        &\geq \frac{\delta_0}{2} \, \lambda_{\max}^{(t)}
    \end{align*}
    Setting $\delta = \delta_0/2$ completes the proof.
\end{proof}

This lemma shows that as the eigenvalues associated to eigenvectors in $\mathcal{N}(S)$ grow, there is at least one off-diagonal $\theta_{ij}$ that grows in absolute value at a comparable rate to $\lambda_{\max}$. This is enough to prove that when $S$ has strictly positive diagonal, a sufficiently strong $\pen$ ensures existence of the penalised likelihood estimate. 
%Note that Lemma \ref{lem2} works with sequences of $\Theta^{(t)}$ as in Lemma \ref{lemma}, without the additional property of Lemma \ref{lemma2}. This will be relevant for results in Section \ref{sec:Combined}.

\begin{myprop}\label{prop:offdiag}
    Let $S$ be PSD with $\dim(\mathcal{N}(S)) = m_0$. Fix $\pen_D = 0$. If $\liminf_{\abs{x} \to \infty} \frac{\pen(x)}{\log\abs{x}} > m_0$  then the penalised likelihood estimate $\hat{\Theta}_{\Pen}$ exists if and only if $S$ has strictly positive diagonal entries.
\end{myprop}
\begin{proof}
    First suppose $S$ has a diagonal entry equal to $0$, w.l.o.g.\ $s_{11} = 0$. Choose $\Theta$ to have eigenvector $v_1 = (1,0,\dots,0)^\top$. Since $v_1 \in \mathcal{N}(S)$, the trace term does not depend on $\lambda_1$ and neither does the penalty term because $v_{1i}v_{1j}=0$ for all $i<j$. So for any fixed $\lambda_2,\dots,\lambda_p$ and $v_2,\dots,v_p$, $l_{\Pen}\left(\,\Theta \mid S \,\right) \to \infty$ as $\lambda_1 \to \infty$.

    Now suppose that all diagonal entries of $S$ are strictly positive. Consider a sequence of $\Theta^{(t)}$ as in Lemma \ref{lemma}. The penalised likelihood is upper bounded by
    $$ l_{\Pen}\left(\,\Theta^{(t)} \mid S \,\right) \leq \sum_{k=1}^{m} \log \lambda_k^{(t)} - \sum_{i < j} \pen\left( \theta_{ij}^{(t)} \right) + c^{(t)}$$
    where $c^{(t)} = \sum_{k=m+1}^{p} \log\lambda_k^{(t)} - \lambda_k^{(t)} v_k^{(t)\top} S v_k^{(t)}$ is uniformly bounded, say $c^{(t)} \leq c$.

    The penalised likelihood can be further upper bounded by removing all penalty terms except that associated to $\theta_{\max}^{(t)}$, since $\pen \geq 0$, and replacing $\lambda_1^{(t)},\dots,\lambda_m^{(t)}$ by $\lambda_{\max}^{(t)}$. To hold uniformly over $V^{(t)}$, this is done when $m=m_0$ is at its maximum value.
    $$l_{\Pen}\left(\,\Theta^{(t)} \mid S \,\right) \leq m_0 \log\lambda_{\max}^{(t)} - \pen\left(\theta_{\max}^{(t)}\right) + c$$
    From Lemma \ref{lem2} we have that $\abs{\theta_{\max}^{(t)}} \geq \delta \, \lambda_{\max}^{(t)}$ for some $\delta>0$ and all $t \geq t_0$. It follows that $\log\abs{\theta_{\max}^{(t)}} \geq \log\lambda_{\max}^{(t)} + \log\delta$. Meanwhile, the condition $\liminf_{\abs{x} \to \infty} \pen(x) / \log\abs{x} > m_0$ implies that there exists $\epsilon>0$ and $x_0$ such that $\pen(x) \geq (m_0+\epsilon)\log\abs{x}$ whenever $\abs{x} > x_0$. These further upper bound the penalised likelihood for large $t$ by
    \begin{align*}
    l_{\Pen}\left(\,\Theta^{(t)} \mid S \,\right) &\leq m_0 \log\lambda_{\max}^{(t)} - (m_0+\epsilon)\log\abs{\theta_{\max}^{(t)}} + c \\
    &\leq m_0 \log\lambda_{\max}^{(t)} - (m_0+\epsilon)\left(\log\lambda_{\max}^{(t)} + \log\delta \right) + c \\
    &= -\epsilon\log\lambda_{\max}^{(t)} -(m_0+\epsilon)\log\delta + c
    \end{align*}
    Since $\lambda_{\max}^{(t)} \to \infty$ and $\delta,c$ are constants, this upper bound tends to $-\infty$.
\end{proof}

When $S$ is a Gaussian sample covariance matrix, the number $m_0$ of eigenvalues equal to 0 is determined with probability 1 by the dimension $p$ and sample size $n$. Additionally, the diagonal entries are positive with probability 1. Hence the condition on $\pen$ ensures existence in the usual continuous data setting.

\begin{mycor}
    If $S$ is a Gaussian sample covariance matrix with $m_0 < p$ eigenvalues equal to 0 and $\liminf_{\abs{x} \to \infty} \frac{\pen(x)}{\log\abs{x}} > m_0$ then the penalised likelihood estimate $\hat{\Theta}_{\Pen}$ exists with probability 1.
\end{mycor}
\begin{proof}
    For unknown $\mu$, $S$ has $m_0 = p-(n-1)$ eigenvalues equal to 0 with probability 1 where $n$ is the sample size. Hence $m_0<p$ implies $n \geq 2$. A diagonal entry of $S$ can be written in terms of $X_1,\ldots,X_n$ as $s_{jj} = \frac{1}{n} \sum_{i=1}^{n} (X_{ij} - \bar{X}_j)^2$ where $X_{ij}$ is the $j$th entry of $X_i$ and $\bar{X}_j = \frac{1}{n} \sum_{i=1}^n X_{ij}$. Hence $s_{jj} = 0$ if and only if $X_{1j} = \cdots = X_{nj}$. Since $X_1,\ldots,X_n$ are independent Gaussian random vectors, this occurs with probability 0 when $n \geq 2$.

    For known $\mu$, instead $S$ has $m_0 = p-n$ eigenvalues equal to 0 with probability 1 so $n \geq 1$. Then $s_{jj} = 0$ if and only if $X_{1j} = \cdots = X_{nj} = \mu$ which occurs with probability 0 when $n \geq 1$.
\end{proof}
The same argument holds for any i.i.d.\ sample from a continuous distribution and so extends to non-Gaussian data.

Proposition \ref{prop:offdiag} shows how to ensure existence for all $S$ with strictly positive diagonal of fixed rank. By taking the smallest possible rank of $S$ (i.e. the largest possible $m_0$), we get a condition that ensures existence for all PSD $S$ with strictly positive diagonal, regardless of rank.

\begin{mycor}\label{cor:offdiagAll}
    Let $S$ be PSD and $pen_D = 0$. If $\liminf_{\abs{x} \to \infty} \frac{\pen(x)}{\log\abs{x}} > p-1$ then the penalised likelihood estimate $\hat{\Theta}_{\Pen}$ exists if and only if $S$ has strictly positive diagonal entries.
\end{mycor}
\begin{proof}
    If $S$ has strictly positive diagonal then it cannot be the zero matrix. It therefore has at most $p-1$ eigenvalues equal to 0. The result follows by using $m_0 = p-1$ in Proposition \ref{prop:offdiag}.
\end{proof}

On the other hand, a weakly growing off-diagonal penalty results in the solution not existing for any only PSD $S$. Similarly to the diagonal penalty, this happens when the penalty is not enough to overpower the log term even when all $p(p-1)/2$ off-diagonals grow.

\begin{myprop}\label{prop:offdiagnonexist}
    Let $q=p(p-1)/2$. If $\pen_D = 0$ and $\liminf_{\abs{x} \to \infty} q \, \pen(x) - \log\abs{x} = -\infty$, then $\hat{\Theta}_{\Pen}$ does not exist for any only PSD $S$.
\end{myprop}
\begin{proof}
    We follow the same style of proof as Proposition \ref{prop:diagnonexist} with $\Theta^{(t)} = I_p + tv_1v_1^\top$ for unit vector $v_1 \in \mathcal{N}(S)$ and $t > 0$. The off-diagonals of $\Theta^{(t)}$ are $\theta_{ij}^{(t)} = t v_{1i} v_{1j}$ and so $\abs{\theta_{ij}^{(t)}} \leq t$. Since $\pen(x)$ is increasing in $\abs{x}$, we have
    \begin{align*}
        \sum_{i<j} \pen\left(\theta_{ij}^{(t)}\right) &\leq \sum_{i<j} \pen(t) \\
        &= q \, \pen(t) \\
        &\leq q \, \pen(1+t)
    \end{align*}
    Hence the penalised likelihood function is lower bounded by 
    $$ l_{\Pen} \left(\,\Theta^{(t)} \mid S \,\right) \geq \log(1+t) - \tr(S) - q \, \pen(1+t)$$
    The condition $\liminf_{x\to\infty} q \, \pen(x) - \log x = -\infty$ ensures that there is an increasing sequence of $t$ such that this lower bound tends to $\infty$ and therefore along this sequence, $l_{\Pen} \left(\,\Theta^{(t)} \mid S \,\right) \to \infty$.
\end{proof}

The results in this section show that when $\pen$ grows super-logarithmically the solution exists for any PSD $S$ with positive diagonal, while if $\pen$ grows sub-logarithmially then the solution does not exist for any only PSD $S$. We have also covered the case when $\pen$ grows at a $c \, \log x$ rate for $c>m_0$ and $c < 1/q$; when $c>m_0$ the solution exists for any $S$ with at most $m_0$ eigenvalues equal to 0 and strictly positive diagonal, while if $c < 1/q$ the solution doesn't exist for any PSD $S$. This leaves the interval $c \in [1/q,m_0]$ where existence is not determined. Existence in this region depends on geometric properties of $\mathcal{N}(S)$. We have assumed worst case scenarios in Proposition \ref{prop:offdiag} where only a single off-diagonal entry grows at the same rate as the maximum eigenvalue and in Proposition \ref{prop:offdiagnonexist} where all off-diagonals grow at the same rate as a single eigenvalue. Under specific assumptions about $\mathcal{N}(S)$, for example on the sparsity of null space vectors, these bounds may be tightened.

\section{Joint diagonal and off-diagonal penalties}\label{sec:Combined}

Section \ref{sec:diag} concentrated on the diagonal penalty while Section \ref{sec:offdiag} focused on the off-diagonal penalty with the other fixed equal to $0$. Since we have assumed non-negative penalty functions, the results of Proposition \ref{prop:diag} and Corollary \ref{cor:DiagProb1} trivially hold for any $\pen \geq 0$ while the existence part of Proposition \ref{prop:offdiag} holds for any $\pen_D \geq 0$, with a sufficiently strong $\pen_D$ also allowing the estimate to exist even when $S$ has diagonal entries equal to $0$ and for lower rank $S$.

However, even when the conditions of both Corollary \ref{cor:DiagProb1} and \ref{prop:offdiag} do not hold, the combination of the diagonal and off-diagonal penalties may still be enough to give existence.

\begin{myprop}\label{prop:joint}
    Let $S$ be PSD with $dim(\mathcal{N}(S)) = m_0$ and strictly positive diagonal entries. Suppose that
\begin{align*}
    \liminf_{x\to\infty}\frac{\mathrm{pen}_D(x)}{\log x} = c_D && \liminf_{\abs{x}\to\infty}\frac{\mathrm{pen}(x)}{\log\abs{x}} = c_O
\end{align*}
If $c_D + c_O > m_0$ then the penalised likelihood estimate $\hat{\Theta}_{\mathrm{Pen}}$ exists.
\end{myprop}
\begin{proof}
Consider a sequence of PD matrices $\Theta^{(t)}$ as in Lemma \ref{lemma}.
%with eigenvalues $\lambda^{(t)}_1,\ldots,\lambda^{(t)}_p$ and eigenvectors $v^{(t)}_1,\ldots,v^{(t)}_p$ ordered such that $v^{(t)}_1,\ldots,v^{(t)}_{m(t)} \in N(S)$ and $v^{(t)}_{m(t)+1},\ldots,v^{(t)}_p \notin N(S)$. Let $\lambda^{(t)}_{\max} = \max_{k=1,\ldots,m(t)}\lambda^{(t)}_k$ be the maximum eigenvalue associated to eigenvectors in $N(S)$, and let $\bar{\lambda}^{(t)}_{\max} = \max_{k=m(t)+1,\ldots,p}\lambda^{(t)}_k$ be the maximum of the remaining eigenvalues. From the discussion in Section 3, we need only consider such sequences with $\lambda^{(t)}_{\max}\to\infty$ and $\sup_t \bar{\lambda}^{(t)}_{\max} < \infty$.
Let $\theta^{(t)}_{\max}$ denote the off-diagonal entry of $\Theta^{(t)}$ with largest absolute value, i.e., $\theta^{(t)}_{\max} = \theta^{(t)}_{i_tj_t}$ where $(i_t,j_t) \in \arg\max_{i < j}\abs{\theta^{(t)}_{ij}}$. By Lemma \ref{lem2}, there exist $\delta > 0$ and $t_0$ such that
$\abs{\theta^{(t)}_{\max}} \geq \delta \lambda^{(t)}_{\max}$
for all $t \geq t_0$.

Next, consider the diagonal entries. Let $v^{(t)}_\max$ be an eigenvector associated with $\lambda^{(t)}_{\max}$. Since $v^{(t)}_\max$ is a unit vector, $\sum_{i=1}^p (v^{(t)}_{\max, i})^2 = 1$, and therefore $\max_{i=1,\ldots,p}(v^{(t)}_{\max, i})^2 \geq 1/p$. Let $i_t = \arg\max_{i=1,\ldots,p}(v^{(t)}_{\max, i})^2$. Then
$$
\theta^{(t)}_{i_ti_t} = \sum_{k=1}^p \lambda^{(t)}_k \left(v^{(t)}_{k i_t}\right)^2 \geq \lambda^{(t)}_{\max}\left(v^{(t)}_{\max, i_t}\right)^2 \geq \frac{\lambda^{(t)}_{\max}}{p}
$$
The penalised likelihood is equal to
$$
l_{\Pen}\left(\,\Theta^{(t)} \mid S \,\right) = \sum_{k=1}^{m} \log\lambda^{(t)}_k - \sum_{i < j}\pen\left(\theta^{(t)}_{ij}\right) - \sum_{i=1}^p\pen_D\left(\theta^{(t)}_{ii}\right) + c^{(t)}
$$
where $c^{(t)} = \sum_{k=m+1}^p \log\lambda^{(t)}_k - \tr(S\Theta^{(t)})$ is uniformly bounded, say $c^{(t)} \leq c$.

Since $\pen$ and $\pen_D$ are non-negative, the penalised likelihood can be upper bounded by removing all penalty terms except those associated to $\theta^{(t)}_{\max}$ and $\theta^{(t)}_{i_ti_t}$, and replacing $\lambda^{(t)}_1,\ldots,\lambda^{(t)}_{m^{(t)}}$ by $\lambda^{(t)}_{\max}$. To hold uniformly over $V^{(t)}$, this is done when $m = m_0$ is at its maximum value. Hence
$$
l_{\Pen}\left(\,\Theta^{(t)} \mid S \,\right) \leq m_0\log\lambda^{(t)}_{\max} - \mathrm{pen}\left(\theta^{(t)}_{\max}\right) - \mathrm{pen}_D\left(\theta^{(t)}_{i_ti_t}\right) + c
$$

The conditions $\liminf_{x\to\infty}\mathrm{pen}_D(x)/\log x = c_D$ and $\liminf_{\abs{x}\to\infty}\mathrm{pen}(x)/\log\abs{x} = c_O$ imply that for any $\epsilon > 0$, there exist thresholds $x_D, x_O$ such that $\pen_D(x) \geq (c_D - \epsilon)\log x$ for all $x \geq x_D$ and $\pen(x) \geq (c_O - \epsilon)\log\abs{x}$ for all $\abs{x} \geq x_O$. Since also $\abs{\theta^{(t)}_{\max}} \geq \delta \lambda^{(t)}_{\max}$ and $\theta^{(t)}_{i_ti_t} \geq \lambda^{(t)}_{\max}/p$, for sufficiently large $t$ we have
\begin{align*}
\mathrm{pen}\left( \theta^{(t)}_{\max} \right) &\geq \left( c_O - \epsilon \right) \log\abs{\theta^{(t)}_{\max}} \geq \left( c_O - \epsilon \right) \left( \log \lambda^{(t)}_{\max} + \log\delta \right) \\
\mathrm{pen}_D\left( \theta^{(t)}_{i_ti_t} \right) &\geq \left( c_D - \epsilon \right) \log\theta^{(t)}_{i_ti_t} \geq \left( c_D - \epsilon \right) \left( \log\lambda^{(t)}_{\max} - \log p \right)
\end{align*}

Substituting these bounds into the upper bound for the penalised likelihood gives
\begin{align*}
l_{\Pen}\left(\,\Theta^{(t)} \mid S \,\right) &\leq m_0\log\lambda^{(t)}_{\max} - \left( c_O - \epsilon \right) \left( \log \lambda^{(t)}_{\max} + \log\delta \right) - \left( c_D - \epsilon \right) \left( \log\lambda^{(t)}_{\max} - \log p \right) + c \\
&= \big(m_0 - (c_O + c_D) + 2\varepsilon\big)\log\lambda^{(t)}_{\max} + (c_O - \varepsilon)\log\delta + (c_D - \varepsilon)\log p + c
\end{align*}

Since $c_D + c_O > m_0$, we can choose $\varepsilon > 0$ small enough such that $m_0 - (c_O + c_D) + 2\varepsilon < 0$. For this choice of $\varepsilon$, the upper bound tends to $-\infty$ as $\lambda^{(t)}_{\max}\to\infty$, which completes the proof.
\end{proof}

This result can be extended to all PSD $S$ by taking the maximum $m_0 = p-1$.

\begin{mycor}
    Let $S$ be PSD with strictly positive diagonal entries. If
    $$
    \liminf_{x\to\infty} \frac{\mathrm{pen}_D(x)}{\log x} + \liminf_{\abs{x}\to\infty} \frac{\mathrm{pen}(x)}{\log\abs{x}} > p - 1
    $$
    then the penalised likelihood estimate $\hat{\Theta}_{\mathrm{Pen}}$ exists.
\end{mycor}
\begin{proof}
If $S$ has strictly positive diagonal then it cannot be the zero matrix. It therefore has at most $p - 1$ eigenvalues equal to $0$. The result follows by applying Proposition \ref{prop:joint} with $m_0 = p - 1$.
\end{proof}

These results show how, even when the individual penalties $\pen$ and $\pen_D$ do not satisfy any of the sufficient conditions for convergence from Sections \ref{sec:diag} and \ref{sec:offdiag}, the combination of the two penalties can still ensure existence. 

Using both penalties similarly affects the case when the solution does not exist for any PSD $S$. The following is an extension of Propositions \ref{prop:diagnonexist} and \ref{prop:offdiagnonexist}.

\begin{myprop}\label{prop:NE_joint}
Let $q = p(p-1)/2$. If
$\liminf_{\abs{x}\to\infty} p\,\pen_D\abs{x} + q \, \pen(x) - \log \abs{x} = -\infty$
then the penalised likelihood estimate $\hat\Theta_{\Pen}$ does not exist for any only PSD $S$.
\end{myprop}
\begin{proof}
We follow the same style of proof as Proposition \ref{prop:diagnonexist} with $\Theta(t) = I_p + tv_1 v_1^\top$ for unit vector $v_1 \in N(S)$ and $t > 0$. Recalling that $\log\det \Theta^{(t)} = \log(1+t)$, $\tr(S\Theta^{(t)}) = \tr(S)$, $\theta_{ii}^{(t)} = 1 + tv_{1i}^2 \leq 1+t$ and $\theta_{ij}^{(t)} = t v_{1i} v_{1j}$ so $\abs{\theta_{ij}^{(t)}} \leq t < 1 + t$, the penalised likelihood function is lower bounded by
\begin{align*}
    l_{\Pen}\left(\,\Theta^{(t)} \mid S \,\right) \geq \log(1+t) - \tr(S) - p \, \pen_D(1+t) - q \, \pen(1+t)
\end{align*}
The condition $\liminf_{x\to\infty} p\,\pen_D(x) + q \, \pen(x) - \log x = -\infty$ ensures that there is an increasing sequence of $t$ such that this lower bound tends to $\infty$ and therefore along this sequence, $l_{\Pen} \left(\,\Theta^{(t)} \mid S \,\right) \to \infty$.
\end{proof}

\section{Generalisations}\label{sec:gen}

In the previous results we made a number of simplifying assumptions about the form of the penalty function. It was assumed to be separable with the same penalty across all diagonal entries and across all off-diagonal entries. These penalties on the individual entries were also assumed to be non-negative, symmetric around 0, increasing in $\abs{\theta_{ij}}$ and continuous, since these are all common properties for penalty functions. However, most of these assumptions were for notational convenience and the results of the previous sections can easily generalise beyond these assumptions.

\subsection{Unequal penalties}\label{subsec:nonequal}

For the more general form of a separable penalty function $\Pen(\Theta) = \sum_{i \leq j} \pen_{ij}(\theta_{ij})$ where a different penalty $\pen_{ij}$ is applied to each $\theta_{ij}$, aside from notational convenience, the results of Sections \ref{sec:diag}-\ref{sec:Combined} easily generalise, although attention must be paid to the precise generalisation of the conditions. In particular, it must be ensured that any bounds in the proof can be maintained under the new condition. 

In existence proofs where an upper bound for the penalised likelihood is required, consider an alternative penalty function $\Pen_\min$ with diagonal penalty $\pen_{\min,D}(x) = \min_{i} \pen_{ii}(x)$ and off-diagonal penalty $\pen_{\min}(x) = \min_{i<j} \pen_{ij}(x)$. Notice that $\Pen_\min$ satisfies the condition of having common penalty functions on the diagonals and on the off-diagonals and that $\Pen_\min \leq \Pen$ so
$$ l_\Pen\left(\,\Theta \mid S \,\right) \leq l_{\Pen_\min}\left(\,\Theta \mid S \,\right)$$
It follows that if $\hat{\Theta}_{\Pen_\min}$ exists then so does $\hat{\Theta}_\Pen$ and so if $\pen_{\min,D}$ or $\pen_\min$ satisfy any of the existence conditions from Sections \ref{sec:diag}-\ref{sec:Combined} then the result extends to $\hat{\Theta}_\Pen$.

On the other hand, the non-existence results require a lower bound on the penalised likelihood which is achieved via an upper bound on the penalty function. For the diagonal penalty this results in a condition on $p \, \pen_D(x)$ while for the off-diagonal penalty there is a condition on $q \, \pen(x)$ where $q$ is the number of off-diagonal entries. For unequal penalties, these are replaced by conditions on the sum of the diagonal penalties $\pen_{\mathrm{sum},D}(x) = \sum_{i=1}^p \pen_{ii}(x)$ and of the off-diagonal penalties $\pen_{\mathrm{sum}}(x) = \sum_{i<j} \pen_{ij}(x)$. The proofs remain unchanged with this replacement.

For simplicity, in the remainder of this section we will revert to the equal penalty $\pen_D$ and $\pen$ framework. However, the results can similarly be generalised.

\subsection{Unbounded below penalties}\label{subsec:nonneg}

The non-negative assumption was again for notational convenience and can easily be relaxed to lower-bounded penalty functions. If the off-diagonal penalty $\pen$ is unbounded below then the solution is guaranteed to not exist because the penalised likelihood function is then unbounded above. Hence it is required that for some $c>-\infty$, $\pen(x) \geq c$ for all $x \in \mathbb{R}$. For simplicity we will continue assuming this holds for $c=0$ so that $\pen$ is non-negative.

However, existence can still be maintained when the diagonal penalty is not lower bounded at $0$, $\lim_{x \to 0^+} \pen_{D}(x) = -\infty$. This is because the log-likelihood also tends to $-\infty$ as any $\theta_{ii} \to 0$. Existence then depends on the limiting behaviour of $\pen_D$ at 0 in comparison to the $\log$ function. Unlike before, existence proofs must now also pay attention to sequences where $\lambda_\min \to 0$. However, this is a necessary but not sufficient condition for $\theta_{ii} \to 0$. If all diagonal entries remain bounded away from 0, then the diagonal penalty remains lower bounded and proving that $l_\Pen\left(\,\Theta \mid S \,\right) \to -\infty$ is as before. Therefore the only additions to the previous proofs are for sequences where some $\theta_{ii} \to 0$. This leads to a direct extension of Proposition \ref{prop:diag}.

\begin{myprop}\label{prop:unbounded}
    If $\lim_{x \to \infty} \pen_D(x) - \log x = \infty$ and $\lim_{x \to 0^+} \pen_D(x) - \log x = \infty$ then $\hat{\Theta}_{\Pen}$ exists for any PSD $S$.
\end{myprop}
\begin{proof}
    As in the proof of Proposition \ref{prop:diag}, removing the trace term and using Hadamard's inequality gives the upper bound
    $$l_{\Pen}\left(\,\Theta \mid S \,\right) \leq \sum_{i=1}^p \log\left(\theta_{ii}\right) - \pen_D\left(\theta_{ii}\right)$$
    From the Proposition conditions we know that if $\theta_{ii} \to 0$ or $\theta_{ii} \to \infty$ then $\log\left(\theta_{ii}\right) - \pen_D\left(\theta_{ii}\right) \to -\infty$. Meanwhile, if $\theta_{ii}$ remains finite then $\log\left(\theta_{ii}\right) - \pen_D\left(\theta_{ii}\right)$ is also finite. It therefore follows that as any $\theta_{ii} \to 0$, the whole sum $\sum_{i=1}^p \log\left(\theta_{ii}\right) - \pen_D\left(\theta_{ii}\right) \to -\infty$, and therefore $l_{\Pen}\left(\,\Theta \mid S \,\right) \to -\infty$ also.
\end{proof}
While this only extends Proposition \ref{prop:diag}, we conjecture that Corollary \ref{cor:DiagProb1} and Proposition \ref{prop:joint} also extend to the case where $\lim_{x \to 0^+} \pen_D(x) - \log x = \infty$ since the same $\theta_{ii} \to 0$ dynamics apply. However, proving these results is slightly more complicated due to dealing with sequences where $\theta_{ii} \to 0$ and $\lambda_\max \to \infty$ together. Such cases require different upper bounds based on the number of $\theta_{ii} \to 0$ and the number of diverging $\lambda_k$.

On the other hand, if $\pen_D$ goes to $-\infty$ at $0$ faster than the $\log$ term, then the penalised likelihood estimate does not exist for any $S$, including PD $S$.
\begin{myprop}
    If $\liminf_{x \to 0^+} \pen_D(x) - \log x = -\infty$ then $\hat{\Theta}_{\Pen}$ does not exist for any PD $S$.
\end{myprop}
\begin{proof}
    Let $\Theta^{(t)} = t I_p$. Then, assuming w.l.o.g.\ that the off-diagonal penalty has $\pen(0)=0$, the penalised likelihood function is
    \begin{align*}
      l_{\Pen}\left(\,\Theta^{(t)} \mid S \,\right) &= p \log t - t \, \tr(S) - p \, \pen_D(t) \\
      &= p \left( \log t - \pen_D(t) \right) - t \, \tr(S)
    \end{align*}
    The condition $\liminf_{x \to 0^+} \pen_D(x) - \log x = -\infty$ ensures that there is a sequence of $t \to 0^+$ such that $\log t - \pen_D(t) \to \infty$. Along this sequence we also have $t \, \tr(S) \to 0$ and so $l_{\Pen} \left(\,\Theta^{(t)} \mid S \,\right) \to \infty$.
\end{proof}

\subsection{Non-symmetric penalties}\label{subsec:nonsym}

Symmetry around $0$ of the off-diagonal penalty function $\pen$ was also for notational convenience. All existence results in Sections \ref{sec:offdiag} and \ref{sec:Combined} extend to the non-symmetric case by lower bounding $\pen$ with $\pen_\min(x) = \min\{ \pen(x), \pen(-x) \}$. The non-existence results also extend with the same replacement, but a slight adaptation of the proof is required to ensure that the chosen sequence of $\Theta^{(t)}$ has diverging $\theta_{ij}^{(t)}$ of the correct sign.

\subsection{Not increasing penalties}\label{subsec:notincreasing}

In penalised likelihoods it is most common to use penalty functions that are increasing in $\abs{\theta_{ij}}$. This is because penalty functions are generally used to regularise and encourage sparsity in the estimate, which is achieved by applying a greater penalty to larger values. However, this is not required for any of the existence proofs and so all existence results still apply to penalty functions that are not increasing, as long as both $\pen_D$ and $\pen$ satisfy the given limiting conditions.

\subsection{Non-continuous penalties}

Most commonly used penalty functions are continuous due to the computational challenges of working with non-continuous penalties. While the proofs of Propositions \ref{prop:diag} and \ref{prop:offdiag} remain valid for non-continuous penalty functions, in so much as they show that the penalised likelihood function tends to negative infinity as the eigenvalues of $\Theta$ grow, this is not enough to prove the existence of the penalised likelihood estimate. This is because the maximum may not be achieved due to the discontinuity.

One special case is the $L_0$ penalty on off-diagonal entries where $\pen(x)=0$ for $x=0$ and $\pen(x)=\rho>0$ for $x \neq 0$. This directly penalises the number of non-zero off-diagonals in $\Theta$. The maximised $L_0$ penalised likelihood under a specific graphical model simply returns the MLE under that model. Hence the $L_0$ penalised likelihood estimate can be seen as comparing the MLEs over all possible graphical models, subtracting a multiple of the number of edges from the maximised likelihood. Since the MLE is guaranteed to exist for all graphical models when $S$ is PD \citep{uhler2012geometry}, this implies the existence of the $L_0$ penalised likelihood estimate for any PD $S$. While a formal proof goes beyond the scope of this paper, we conjecture that the $L_0$ penalised likelihood estimate will continue to exist for positive semi-definite $S$ when combined with a suitably strong penalty on the diagonal entries, satisfying the condition of Proposition \ref{prop:diag}.

\subsection{Partial correlation based penalties}\label{subsec:PC}

Recently a new class of penalised likelihoods for Gaussian graphical models have been introduced that directly penalise the partial correlations rather than the off-diagonal entries of $\Theta$ \citep{Carter2023,bogdan2025identifying,carter2025existence}. We denote the partial correlations by $\Delta_{ij} = -\theta_{ij}/\sqrt{\theta_{ii}\theta_{jj}}$ which must be between $-1$ and $1$ for the resulting matrix to be PD. A penalty function is then placed on $\Delta = (\Delta_{ij})$ and the diagonal entries of $\Theta$, $\theta = \mathrm{diag}(\Theta)$
$$\Pen(\Delta,\theta) = \sum_{i=1}^p \pen_D(\theta_{ii}) + \sum_{i<j} \pen(\Delta_{ij})$$
The results of Section \ref{sec:diag} directly generalise to this class of penalised likelihoods, as long as $\pen(\Delta_{ij})$ satisfies the usual conditions of being lower bounded and continuous.

So far only an $L_1$ penalty has been implemented for the partial correlations, which is bounded on the interval $[-1,1]$. For any such bounded penalty, existence of the solution relies solely on the diagonal penalty $\pen_D$. However, even if the penalty on $\Delta_{ij}$ is allowed to grow infinitely such that $\pen(\Delta_{ij}) \to \infty$ as $\abs{\Delta_{ij}} \to 1$, it is likely that this still isn't enough to ensure existence. This is because there are only PSD $\Delta$ such that no $\abs{\Delta_{ij}} = 1$ and so $\abs{\Delta_{ij}} \to 1$ does not correspond to the boundary of PD matrices. It therefore seems that partial correlation based penalties will always rely on penalisation of the diagonal or other similar quantities in order to exist for only PSD $S$.

\section{Bayesian methods and posterior propriety}\label{sec:Bayes}

For a prior distribution $p(\Theta)$, the posterior given $S$ from a sample of size $n$ is
$$ p\left(\,\Theta \mid S \,\right) \propto \det(\Theta)^{n/2} \exp\!\left(-\frac{n}{2}\,\tr(S\Theta)\right) \, p(\Theta) \, \mathbb{I}(\Theta \succ 0)$$
By taking the logarithm and considering the penalty function $\Pen(\Theta) = -\frac{2}{n} \log p(\Theta)$, the posterior can be related to the penalised likelihood function via
$$ \log p\left(\,\Theta \mid S \,\right) = \frac{n}{2} \, l_\Pen\left(\,\Theta \mid S \,\right) + c$$
for some constant $c$. It follows that the maximum a posteriori (MAP) estimate under $p(\Theta)$ is equal to the penalised likelihood estimate with penalty $\Pen$. The results of the previous sections can therefore be directly applied to prove the existence of the MAP estimate under separable priors of the form
$$ p(\Theta) \propto \prod_{i=1}^p \pi_D(\theta_{ii}) \prod_{i < j} \pi(\theta_{ij}) \, \mathbb{I}(\Theta \succ 0) $$
by substituting $\pen_D(x) = -\frac{2}{n} \log \pi_D(x)$ and $\pen(x) = -\frac{2}{n} \log \pi(x)$, noting that the assumptions of $\pen_D$ and $\pen$ should be passed on to $\pi_D$ and $\pi$ respectively. In this section we will continue to assume that $\pi_D$ and $\pi$ are continuous.

Beyond existence of the MAP, it is useful in Bayesian statistics to ensure that the posterior distribution is proper, i.e. the unnormalised posterior has finite integral. The existence proofs for the penalised likelihood estimate demonstrated that, under certain conditions on $\Pen$, $l_\Pen\left(\,\Theta \mid S \,\right) \to -\infty$ as $\lVert \Theta \rVert \to \infty$. It therefore follows that the posterior $p\left(\,\Theta \mid S \,\right) \to 0$ as $\lVert \Theta \rVert \to \infty$ under the equivalent conditions on $p(\Theta)$. However, this is not enough to ensure that the posterior is proper because the integral of a function over an unbounded set can still be infinite even if the function is finite and tends to 0 at the boundary.

To investigate the propriety of the posterior further, we begin with the case where $S$ is PD.

\begin{myprop}\label{prop:ProperPD}
    Suppose $S$ is PD, $n \geq 1$, $\pi$ is bounded, and $\pi_D$ has $\limsup_{x \to \infty} \pi_D(x) < \infty$ and $\int_{0}^1 x^{(n+p-1)/2} \pi_D(x) \, dx < \infty$. Then the posterior $p\left(\,\Theta \mid S \,\right)$ is proper.
\end{myprop}
\begin{proof}
    We want to show that $Z < \infty$ where
    $$ Z = \int_{\Theta \succ 0} \det(\Theta)^{n/2} \exp\!\left(-\frac{n}{2}\,\tr(S\Theta)\right) \, \prod_{i=1}^p \pi_D(\theta_{ii}) \prod_{i < j} \pi(\theta_{ij}) \, d\Theta$$
    Using Hadamard's inequality, $\det(\Theta)^{n/2} \leq \prod_{i=1}^p \theta_{ii}^{n/2}$. Since $S$ is PD, the trace term is lower bounded by $\tr(S\Theta) \geq \lambda_\min (S)\tr(\Theta) = \lambda_\min(S) \sum_{i=1}^p \theta_{ii}$. Since $\pi$ is bounded, $\prod_{i < j} \pi(\theta_{ij}) < M$ for some constant $M<\infty$. Hence we have the upper bound
    $$ Z \leq \int_{\Theta \succ 0} \prod_{i=1}^p \theta_{ii}^{n/2} \exp\!\left(-\frac{n}{2}\,\lambda_\min(S) \sum_{i=1}^p \theta_{ii}\right) \prod_{i=1}^p \pi_D(\theta_{ii}) \, M \, d\Theta$$
    Note that the integrand does not depend on the off-diagonals $\theta_{ij}$ and for given diagonals $\theta=(\theta_{ii})_{i=1}^p$, the off-diagonals are restricted by positive definiteness to $\abs{\theta_{ij}} < \sqrt{\theta_{ii}\theta_{jj}}$. For fixed diagonals, the volume of the space of off-diagonals is therefore upper bounded by $\prod_{i<j} 2\sqrt{\theta_{ii}\theta_{jj}} = 2^q \prod_{i=1}^p \theta_{ii}^{(p-1)/2}$ where $q = p(p-1)/2$ is the number of off-diagonals. Hence
    \begin{align*}
        Z &\leq C \int_{\theta \in (0,\infty)^p} \prod_{i=1}^p \theta_{ii}^{(n+p-1)/2} \exp\!\left(-\frac{n}{2}\,\lambda_\min \sum_{i=1}^p \theta_{ii}\right) \prod_{i=1}^p \pi_D(\theta_{ii}) \, d\theta \\
        &= C \prod_{i=1}^p \int_0^\infty x^{(n+p-1)/2} \exp(-ax) \pi_D(x) \, dx
    \end{align*}
    where $C=2^q M$ and $a = \frac{n}{2}\lambda_\min$ are constants.

    On $[1,\infty)$, $\pi_D$ is bounded (by continuity and $\limsup_{x \to \infty} \pi_D(x) < \infty$) and $\int_1^\infty x^{(n+p-1)/2} \exp(-ax) \, dx < \infty$. On $(0,1]$, $\exp(-ax) \leq 1$ and $\int_0^1 x^{(n+p-1)/2} \pi_D(x) \, dx < \infty$. Hence $Z < \infty$ and so the posterior is proper.
\end{proof}

In terms of the corresponding penalty functions, the conditions basically mean that $\pen_D(x)$ and $\pen(x)$ cannot tend to $-\infty$ as $\abs{x} \to \infty$, while close to $0$, $\pen_D$ cannot go to $-\infty$ faster than $c \log x$ where $c = 1 + (p+1)/n$. In other words, $\limsup_{x \to 0^+} \pen_D(x) / \log x < c$. Note that this is actually a weaker condition than in Section \ref{subsec:nonneg} where we required $\lim_{x \to 0^+} \pen_D(x) - \log x = \infty$. This is because a function can tend to $\infty$ at 0 (i.e. the maximum doesn't exist) but still have finite integral on $(0,1)$.

Moving on to PSD $S$, we encounter the same problem as with penalised likelihoods where the trace term can remain constant when $\Theta$ grows in directions in the null space of $S$. This requires stronger conditions on $\pi$ and $\pi_D$ to maintain propriety of the posterior for all PSD $S$ and any $n$.

\begin{myprop}\label{prop:ProperPSD}
    Let $S$ be PSD and $n\geq1$. If $\pi$ is proper, $\int_{-\infty}^\infty \pi(x) \, dx < \infty$, and $\pi_D$ has finite $n/2$ moment, $\int_0^\infty x^{n/2}\pi_D(x) \, dx < \infty$, then the posterior $p\left(\,\Theta \mid S \,\right)$ is proper.
\end{myprop}
\begin{proof}
    Let $Z$ be as in the previous proof. Using the same Hadamard inequality and $\exp\left(-\frac{n}{2}\tr(S\Theta)\right) \leq 1$ we have
    \begin{align*}
        Z &\leq \int_{\Theta \succ 0} \prod_{i=1}^p \theta_{ii}^{n/2} \pi_D(\theta_{ii}) \prod_{i<j} \pi(\theta_{ij}) \, d\Theta \\
        &\leq \prod_{i=1}^p \int_0^\infty \theta_{ii}^{n/2} \pi_D(\theta_{ii}) \, d\theta_{ii} \prod_{i<j} \int_{-\infty}^\infty \pi(\theta_{ij}) \, d\theta_{ij}
    \end{align*}
    Each of these integrals is finite by the conditions on $\pi_D$ and $\pi$. Hence $Z < \infty$ and so the posterior is proper.
\end{proof}
Here we used a finite $n/2$ moment for $\pi_D$ which implies that $\int_0^1 x^{n/2} \pi_D(x) \, dx < \infty$ also. However, this is a slightly stronger condition than that in Proposition \ref{prop:ProperPD} where $\int_0^1 x^{(n+p-1)/2} \pi_D(x) \, dx < \infty$. This is because we haven't used the upper bounded condition on $\pi$, relying only on propriety. But $\pi$ being proper is actually a weaker condition around $0$ than boundedness, because a proper $\pi$ can still have $\lim_{x \to 0} \pi(x) = \infty$. It therefore seems likely that this weaker condition on $\pi_D$ is also enough for the PSD case in Proposition \ref{prop:ProperPSD}, since the trace term, and therefore $S$, does not affect the integral around 0. 

In terms of the corresponding penalty functions, propriety of $\pi$ is given by $\liminf_{\abs{x} \to \infty} \pen(x) / \log x > 2/n$, while $\pi_D$ having finite $n/2$ moment is achieved if $\limsup_{x \to 0^+} \pen_D(x) / \log x < 1 + 2/n$ and $\liminf_{x \to \infty} \pen_D(x) / \log x > 1 + 2/n$. In particular, this $x \to \infty$ condition on $\pen_D$ is stronger than that of Proposition \ref{prop:diag}.

\section*{Acknowledgements}
This research was supported by the EUTOPIA Science and Innovation Fellowship Programme and funded by the European Union Horizon 2020 programme under the Marie Skłodowska-Curie grant agreement No 945380.

\bibliographystyle{plainnat}
\bibliography{myrefs}

\end{document}